\documentclass[12pt]{article}
\usepackage{amsmath, amsthm, amscd, amssymb}

\theoremstyle{plain}

\theoremstyle{definition}

\theoremstyle{remark}

\newcommand{\N}{\mathbb{N}}
\newcommand{\Z}{\mathbb{Z}}

\newcommand{\R}{\mathbb{R}}

\newcommand{\ba}{\begin{aligned}}
\newcommand{\ea}{\end{aligned}}
\newcommand{\bt}{\begin{thm}}
\newcommand{\et}{\end{thm}}
\newcommand{\bc}{\begin{corollary}}
\newcommand{\ec}{\end{corollary}}
\newcommand{\bl}{\begin{lemma}}
\newcommand{\el}{\end{lemma}}
\newcommand{\bpf}{\begin{proof}}
\newcommand{\epf}{\end{proof}}
\newcommand{\bpb}{\begin{problem}}
\newcommand{\epb}{\end{problem}}
\newcommand{\bd}{\begin{definition}}
\newcommand{\ed}{\end{definition}}
\newcommand{\bn}{\begin{note}}
\newcommand{\en}{\end{note}}
\newcommand{\bp}{\begin{proposition}}
\newcommand{\ep}{\end{proposition}}
\newcommand{\be}{\begin{example}}
\newcommand{\ee}{\end{example}}
\newcommand{\bex}{\begin{exercise}}
\newcommand{\eex}{\end{exercise}}
\theoremstyle{plain}
\newtheorem{thm}{Theorem}[section]

\newtheorem{lemma}[thm]{Lemma}
\newtheorem{corollary}[thm]{Corollary}
\newtheorem{proposition}[thm]{Proposition}
\newtheorem{exercise}[thm]{Exercise}
\newtheorem{problem}[thm]{Problem}

\theoremstyle{definition}
\newtheorem{definition}[thm]{Definition}
\newtheorem{remark}[thm]{Remark}
\newtheorem{example}[thm]{Example}

\newtheorem {note}{Note}[section]
\swapnumbers
\theoremstyle{remark}

\theoremstyle{plain}

\title{Discrete/Continuous Elliptic Harnack Inequality and Kernel Estimates for Functions of the Laplacian on a Graph}
\author{Mark Cerenzia and Laurent Saloff-Coste}
\date{\today}

\begin{document}
\maketitle

\begin{abstract}
This paper introduces certain elliptic Harnack inequalities for harmonic functions in the setting of the product space $M \times X$, where $M$ is a (weighted) Riemannian Manifold and $X$ is a countable graph. Since some standard arguments for the elliptic case fail in this ``mixed" setting, we adapt ideas from the discrete parabolic case found in the paper \cite{Delmotte3} of Delmotte. We then present some useful applications of this inequality, namely, a kernel estimate for functions of the Laplacian on a graph that are in the spirit of the paper \cite{CGT1} by Cheeger, Gromov, and Taylor. This application in turn provides sharp estimates for certain Markov kernels on graphs, as suggested in Section $4$ of the paper \cite{PDSC1} by Persi Diaconis and the second author. We then close with an application to convolution power estimates on finitely generated groups of polynomial growth.
\end{abstract}

In \cite{Moser1}, Moser proves a (by now classical) elliptic Harnack inequality for positive solutions of a uniformly elliptic differential operator. In its most basic form, the elliptic Harnack inequality in $\R^m$ says that, for $\delta \in (0,1)$, there is a $C= C(\delta) >0$ such that if $f$ is a positive harmonic function on $B = B_r(x) \subset \R^m$, then
$$
\sup_{\delta B} f \leq C \min_{\delta B} f
$$
Here, $\delta B$ denotes the ball concentric with $B$ and with radius $\delta r$. Moser later extends this result in \cite{Moser2} to derive the more powerful parabolic Harnack inequality. Inequalities of this type have many applications, such as exhibiting Holder continuity, deriving Gaussian bounds for the heat kernel, bounding the dimension of the space of harmonic functions, as well as many other results concerning the underlying geometry of spaces. 

It was shown by the second author in \cite{SC2} (and independently by Grigor'yan \cite{Grigor1}) that the conjunction of two geometric properties, Volume Doubling and the $L^2$ Poincar\'{e} Inequality, is equivalent to the parabolic Harnack inequality. These concepts in turn are equivalent to the heat kernel satisfying Gaussian bounds. This circle of ideas was later extended to the graph setting by Thierry Delmotte, first for the elliptic case in \cite{Delmotte1} and then later for the parabolic case in \cite{Delmotte3}.

Here, we are interested in the analogous elliptic Harnack inequality for harmonic functions on a product space $M \times X$, where $M$ is a (weighted) Riemannian manifold and $X$ is a locally finite graph (see Section \ref{Notation and Set Up} for definitions). As we will see below, this result is useful for extending traditional applications of Harnack inequalities to discrete spaces (see Sections \ref{Kernel Estimates on Graphs} and \ref{applications}). Although a good portion of the standard arguments of Moser and Delmotte can be used in deriving the elliptic Harnack inequality for our ``mixed" space, we will see below that the arguments do not extend as readily as one may initially expect. In fact, a crucial argument for the discrete elliptic case \cite{Delmotte1} does not work in this mixed setting.

The outline of this paper is as follows. The first section will introduce the notation and relevant geometric assumptions for our main results. Here we discuss classical geometric assumptions as well as an additional uniformity assumption required for the graph setting. In section \ref{ellipticsection}, we address the issues raised in the previous paragraph to confirm that the elliptic Harnack inequality does indeed extend to the mixed product $M \times X$. Presenting our main application, Section \ref{Kernel Estimates on Graphs} derives kernel estimates for functions of the Laplacian on a graph. In pursuing this application, we sharpen an interesting and amazingly general bound of Carne and Varopoulos in \cite{Carne1}. Lastly, we mention some corollaries of these results for Markov Chains on Groups.  

\section{Notation and Set Up}
\label{Notation and Set Up}

Part of this paper concerns estimates on the product of a (Riemannian) manifold and a (symmetrically) weighted graph. This section covers the main definitions related to a symmetrically weighted graph, which is akin to a Riemannian structure and which naturally gives rise to a reversible Markov chain (see Section \ref{Kernel Estimates on Graphs}).

\bn
It is important to emphasize that as we make estimates, a given constant will absorb extraneous factors. Although some computations with such constants are made explicit (e.g., by including powers of the constant after a step), it is implied that these will be absorbed into a single constant, usually having the same name.
\en

The continuous part of our product will be a complete Riemannian Manifold $(M, g)$. We denote the resulting volume measure by $\mu$ and define the Laplacian on $C_c^{\infty}(M)$ as $\Delta = \text{div grad }.$ This operator extends to an (essentially) self-adjoint operator on $L^2(M)$. The reader should take care to note our sign convention, which agrees with the discrete Laplacian we define below. In particular, $\Delta$ defined above is a \emph{nonpositive} operator. In addition, our results will apply to weighted spaces. That is, a space with volume measure $d \nu = \sigma d\mu$ for some positive function $\sigma \in C^{\infty}(M)$. The associated Laplacian then becomes
$$
\Delta_{\sigma} = \sigma^{-1} \text{div} ( \sigma \text{ grad } )
$$
(see, for example, the paper of Grigor'yan and Saloff-Coste \cite{Grigor2}). 

For the other component of the product, let $X$ be a discrete space and let $\mu : X \times X \to \R^+$, $(u,v) \mapsto \mu_{uv}=\mu_{vu}$ be a symmetric weight. This function induces a graph structure on $X$ by declaring an edge between $u$ and $v$, written $u \sim v$, if $\mu_{uv} >0$. Define a measure by $m(u) = \sum_{u \sim v} \mu_{uv}$ and the ``volume" of subsets $E \subset X$ by
$$
m(E) : = \sum_{u \in E} m(u).
$$
We can further define the (geodesic) distance between two vertices $u, v \in X$ as $d_X(u,v) = $ shortest number of edges between $u$ and $v$. We can then denote by $B_r(u)$ the \emph{open} ball $\{ v \in X | \ d_X(u,v) < r\}$. Notice that our convention differs from the notation of Delmotte, who lets $B(u,r)$ denote the \emph{closed} ball $\{ v \in X | d_X(u,v) \leq r\}$.
For simplicity, we assume the graph is connected so that $m(u) > 0$ for all $u \in X$. Lastly, for a subset $E \subset X$, define its boundary as $\partial E = \{ u \notin E   \ | u \sim v, \text{for some} \ v \in E \}$ and set $E^* = E \cup \partial E$.

As the name of this paper suggests, we aim to discuss analytic results that involve this graph setting. Thus we recall a few relevant quantities of the discrete calculus. Let $h,g : X \to \R$ be functions on our graph. For the gradient, we will write
$$
\nabla_X h \cdot \nabla_X g(u) : =\frac{1}{m(u)} \sum_{v \in X} \mu_{uv}(h(u)-h(v))(g(u)-g(v)).
$$
and for the discrete Laplacian, 
$$
\Delta_{X} h(u) : = \frac{1}{m(u)} \sum_{v \sim u} \mu_{uv} (h(v) - h(u))
$$
Notice these quantities are well defined since the graph is assumed to be connected. 

Consider the mixed space $M \times X$ endowed with the Pythagorean metric $d:= \sqrt{d_M^2 + d_{X}^2}$ (this choice is only for convenience). For a function $f :M \times X \to \R$, $(x,u) \mapsto f(x,u)$, we write
$$
|\nabla f(x,u)|^2 = |\nabla_M f(x,u)|^2+|\nabla_{X} f(x,u)|^2.
$$
Likewise, the \emph{Laplacian} on $M \times X$ becomes
$$
\Delta f(x,u) : = \Delta_M f(x,u) +\Delta_X f(x,u).
$$
Here, the operators $\nabla_M, \Delta_M, \nabla_X, \Delta_X$ are applied to the indicated argument; for example,
$$
\Delta_X f(x,u) = \sum_{v \sim u} \frac{\mu_{uv}}{m(u)} (f(x,v)-f(x,u)).
$$

We close this section with a discussion of harmonic functions in our setting. Let $\Omega = \Omega_M \times \Omega_X \subset M \times X$ be a subset with $\Omega_M$ open. We now define the appropriate notion of solution to the equation
$$
\Delta f(x,u) = \Delta_M f(x,u)+\Delta_X f(x,u)= 0
$$
for $(x,u) \in \Omega$. 
\bd \label{harmonic}
We say that $f: \Omega = \Omega_M \times \Omega_X \to \R$ is $\mathit{(weakly) \ harmonic}$ on $\Omega$ if $f_u \in L^2(\Omega_M, \R)$ for every $u \in \partial \Omega_X$, $f_u \in H^1(\Omega_M, \R)$ for every $u \in \Omega_X$, and
$$
 \sum_{u \in \Omega_X} \int_{\Omega_M}  \left [ - \nabla_M f_u(x) \cdot \nabla_M \phi_u(x) + \phi_u(x)\Delta_X f_u(x)    \right ] = 0
$$
for all $\phi \in H^1_0(\Omega, \R)$.
\ed
Fix $u \in \Omega_X$, and write $f_u(x):= f(x,u)$ and $g_u(x) :=  - \sum_{v \sim u} \frac{\mu_{uv}}{m(u)} f_v(x)$. Define a second order elliptic operator $L = \Delta_M - I$, where $I$ is the identity operator. Then the equation above becomes
$$
(Lf_u)(x) = g_u(x).
$$
With this notation, the regularity of a given harmonic function $f$ on $\Omega$ can easily be determined by the classical elliptic regularity theorem (see, e.g., Chapter 9 of Folland's text \cite{Folland}). Namely, if $Lf_u = g_u \in H^s(\Omega_M, \R)$, $s \geq 0$, and if $f_u \in H^1(\Omega_M, \R)$, then this theorem guarantees that $f_u \in  H^{s+2}(\Omega_M,\R)$. Thus, if $\Omega_X = X$, then the solution $f$ has enough room to completely ``self-improve" until $f_u \in  C^{\infty}(\Omega_M,\R)$ for all $u \in \Omega_X$. If, on the other hand, $\Omega_X$ is a proper subset of $X$, then the regularity $f_u$ depends both on the distance of $u$ from the boundary $\partial \Omega_X$ and on the regularity of $f_v$ for $v \in \partial \Omega_X$. To see this, note that the regularity of $g_u$ depends on the regularity of $f$ at neighbors $v$ of $u$. Since the smallest of these degrees of regularity determines the regularity of $g_u$, it also determines the regularity of $f_u$. We summarize the above argument precisely in the following proposition.
\bp
Fix $s \geq 0$. Assume that $f$ is weakly harmonic on $\Omega$ and that $f_u \in H^s(\Omega_M, \R)$ for $u \in \partial \Omega_X$. Fix $u_0 \in \Omega_X$ and let $p = d_X(u_0, \partial \Omega_X)$, where we allow the possible value $\infty$. Then $f_{u_0} \in  H^{s+2p}(\Omega_M,\R)$. 
\ep
Lastly, a note on existence. If $\Omega_X$ is finite, then we can dually view our equation
$$
\Delta_M f_u(x) - f_u(x) + \sum_{v \sim u} \frac{\mu_{uv}}{m(u)} f_v(x) = 0
$$
as a nonlinear elliptic system of differential equations. The monograph \cite{ChenWu1} by Chen and Wu is a good reference for the theory of such systems. Indeed, Chapter 11 covers an existence theorem based on variational arguments that works well for our setting. In addition, the reader should note that the proof of Theorem \ref{kernel} below contains a construction of a harmonic function on a mixed space with an interval for its continuous component.

\section{Geometric Properties}
\label{Geometric Properties}

Now that the notation is established, we may begin discussing the geometric properties essential to the elliptic Harnack inequality. We say a Riemannian manifold $(M, g)$ satisfies the \emph{Volume Doubling} property $V(D,r_0)$ if
$$
\exists D : \ \mu( B(x,2r)) \leq D \mu (B(x,r))
$$ 
and the (weak) Poincar\'{e} Inequality $P(P, r_0)$ if
$$
\exists P : \ \int_{B(x,r)} |f-f_{x,r}|^2 dV \leq P r^2 \int_{B(x,2r)} |\nabla f |^2 dV,
$$
where in each case $0 < r<r_0$ and $f$ is any smooth function on $B(x,2r)$. In the manifold setting, it is a fact (see \cite{SC1}, Corollary 5.3.5) that the conjunction of volume doubling and the weak Poincar\'{e} Inequality is equivalent to the strong Poincar\'{e} Inequality (where $2B$ is replaced by $B$ in the integrand on the right).

The above two definitions are the same on a graph, but our discrete component must satisfy an additional property. We say a weighted graph $(X, \mu_{uv}) $ satisfies $\Delta^*(\alpha)$ if
$$
u \sim v \implies \mu_{uv} \geq \alpha m(u)
$$
This property implies that the graph is \emph{locally uniformly finite}, i.e., there exists a constant $A>0$ such that $|\{ v : v \sim u \} |\leq A$, for all $u \in X$. In particular, we may take $A = \alpha^{-1}$. (It is worth noting that if the graph satisfies $V(D_X, r_0)$, then it is locally uniformly finite with $A = D_X^2$.) In \cite{Delmotte2}, Delmotte discusses the extra assumption of requiring loops at each vertex, i.e., $u \sim u$. He lets $\Delta(\alpha)$ denote the conjunction of requiring such loops and the property $\Delta^*(\alpha)$. Together, these assumptions allow a comparison of the discrete kernel with a constructed continuous kernel (see \cite{Delmotte2, Delmotte3} for details). Although we use ideas found in the manipulation of this continuous kernel, we will not need to perform such a comparison, and as a result requiring loops is unnecessary for us.

An elliptic Harnack inequality holds on a Manifold that satisfies Volume Doubling and the (weak) Poincar\'{e} Inequality \cite{SC1}. Thierry Delmotte proves an elliptic Harnack inequality in \cite{Delmotte1} on a weighted graph with all three of the above properties. In this paper, we wish to address whether the product of two such spaces also supports an elliptic Harnack inequality. 

The reader may wonder at first why it is not sufficient to prove that the product space inherits the geometric properties of its components. Perhaps then one can derive the desired inequality just as for each individual component. Unfortunately, as we will see, the elliptic Harnack does not follow so readily. Nevertheless, we still intend to use this inheritance of geometry, whose confirmation is contained in the two lemmas below.

As discussed above, let $(M, g)$ be a Riemannian manifold with induced measure $\mu$ and let $(X, \mu_{uv})$ be a symmetrically weighted graph with induced measure $m$. Fix $r_0 \in (0,\infty]$ (in many applications, $r_0$ will in fact be $\infty$). Assume that $M$ and $X$ satisfy the Volume Doubling Property, $V(D_M, r_0)$ and $V(D_X, r_0)$, respectively, as well as the Poincar\'{e} Inequality, $P(P_M, r_0)$ and $P(P_X, r_0)$, respectively. Let $\pi$ denote the product measure on $\Pi = M \times X$, i.e., $d\pi = d\mu \times  dm$, and let $d$ be the Pythagorean product distance $d:= \sqrt{d_M^2 + d_{X}^2}$. We will start by showing that the space $(\Pi, d, \pi)$ inherits these two crucial properties from its components.
\bl \label{vdlemma}
The product $\Pi$ satisfies  $V(D, r_0)$, where
$$
D=D_M D_X \left ( 2 \sqrt{2} \right )^{\frac{\log D_M+ \log D_X}{\log 2}}.
$$
\el
\bpf
It is a straightforward computation to show that volume doubling implies, for $x \in M$ and $r \geq s$,
\begin{equation} \label{vdeqn}
\mu(B(x,r)) \leq D_M\left ( \frac{r}{s}  \right )^{\log D_M/\log 2} \mu(B(x,s)),
\end{equation}
and similarly for $X$. As in the statement, let
$$
V=D_M D_X \left ( 2 \sqrt{2} \right )^{\frac{\log D_M+ \log D_X}{\log 2}}.
$$
Then with $2r \geq \frac{r}{\sqrt{2}}$, we compute, for $(x,y) \in \Pi$,
$$
\ba
\pi(B((x,y), 2r) & \leq \mu(B(x, 2r))m(B(y, 2r)) \\ &
\leq V \mu(B((x, r/\sqrt{2}))m(B(y, r/\sqrt{2})) 
\\ & \leq V \pi(B((x,y), r)).
\ea
$$
For the first and last inequalities, we have used the uniform equivalence of the natural product metric $d$ on $\Pi$ and the max metric $\rho := \max \{ d_M, d_X\}$:
$$
\frac{1}{\sqrt{2}} d \leq \rho \leq d.
$$
\epf
Analogously, the Poincar\'{e} Inequality also holds on our the product space.
\bl
The product $\Pi= M \times X $ satisfies $P(P, r_0)$ with 
$$
P = 2 P_M P_X.
$$
\el
\bpf
For the first step, consider $B=B_M \times B_M \subset \Pi$, where the component balls have radius $0<r<r_0$. Remember that both $M$ and $X$ satisfy the strong Poincar\'{e} Inequality (\cite{SC1}, Corollary 5.3.5). Now let $f_B$ denote the the average of $f$ over $B$:
$$
f_B:= \frac{1}{\pi(B)} \int_B f d \pi.
$$
Notice that by Fubini's Theorem, $f_B = (f_{B_M}(u))_{B_X}$, where the averages are taken with respect to the obvious component. Note also that by Jensen's Inequality, $\phi(f_{B_M}(u)) \leq (\phi \circ f(u))_{B_M}$ for every convex $\phi: \R \to \R$. Recall, lastly, the inequality $|a-c|^2 \leq 2( |a-b|^2+|b-c|^2)$. Keeping these facts in mind, we compute

\begin{eqnarray*} 
 \int_{B_X} \int_{B_M} && |f(x,u)-f_B|^2 d\mu(x) dm(u)  \leq  \\ 
&& 2 \int_{B_X} \int_{B_M} (|f(x,u)-f_{B_M}(u)|^2 + |f_{B_M}(u)-f_B|^2) d\mu(x) dm(u). 
\end{eqnarray*} 
To bound the first term of the integrand on the right, we invoke $P(P_M,r_0)$ in $N$:
$$
\int_{B_M} |f(x,u)-f_{B_M}(u)|^2 d\mu(x) \leq P_M r^2 \int_{B_M} | \nabla_M f(x,u)|^2 d\mu(x)
$$
For the other term, we use $P(P_X,r_0)$ in $X$ and Jensen's with $| \cdot |^2$:

\begin{equation}
\begin{split} \nonumber
\int_{B_X} \int_{B_M} |f_{B_M}(u)-f_B|^2 & d\mu(x) dm(u)  \\
&= \mu(B_M) \int_{B_X} |f_{B_M}(u)-(f_{B_M}(u))_{B_X}|^2 dm(u) \\
& \leq \mu(B_M) P_X r^2 \int_{B_X} | \nabla_X f_{B_M}(u)|^2 dm(u) \\
& \leq  P_X r^2 \int_{B_X} \int_{B_M} | \nabla_X f(x,u)|^2 d \mu(x) \ dm(u). \\
\end{split}
\end{equation}

Putting these results together yields
$$
\ba
\int_{B(x,r)} |f(x,u)-f_B|^2 d\pi& \leq \int_{B_X} \int_{B_M} |f(x,u)-f_B|^2 d\mu(x) dm(u) \\ & \leq 2 P_M P_X \int_{B_X} \int_{B_M} | \nabla f|^2 d\mu dm \\
& \leq 2 P_M P_X \int_{B(x, \sqrt{2} r)} | \nabla f|^2 d\pi \\ & \leq 2 P_M P_X \int_{B(x, 2 r)} | \nabla f|^2 d\pi.
\ea
$$
This proves the (weak) Poincar\'{e} Inequality.

\epf

\begin{remark}
As mentioned before, the strong Poincar\'{e} inequality follows from the conjunction of volume doubling and the weak form. The standard proof of this result uses a Whitney Covering argument due to D. Jerison (see Chapter 5 of \cite{SC1} for details). For the graph setting, the proof uses the additional uniformity assumption $\Delta^*(\alpha)$; see \cite{Delmotte3}. See also Section $5.3.2$ of \cite{SC1} for a more general discussion.
\end{remark}

\section{Elliptic Harnack inequality on mixed spaces}
\label{ellipticsection}

For this section, let $M$ be an d-dimensional manifold with volume measure $d\mu$ and $X$ a graph with measure $m(u) = \sum_{v \sim u} \mu_{uv}$, where $\mu_{ij} = \mu_{ji} \geq 0$ is a symmetric weight on $X \times X$. The manifold has the usual geodesic distance $d_M$, and similarly, the graph has the metric $d_{X}(u,v) := $ the smallest number of edges between $u$ and $v$. As above, fix $r_0 \in (0,\infty]$ ($r_0$ will often be $\infty$ in practice). Assume that the components satisfy the Volume Doubling Property, $V(D_M, r_0)$ and $V(D_X, r_0)$, respectively, as well as the Poincar\'{e} Inequality, $P(P_M, r_0)$ and $P(P_X, r_0)$, respectively. As proven in the last section, the product space $\Pi = M \times X$ satisfies $V(D, r_0)$ and $P(P,r_0)$ with respect to the measure $\pi: = \mu \times m$.

We wish to establish the following theorem:

\bt \label{elliptic}
Suppose that $M$ and $X$ are as above and in addition that $X$ satisfies $\Delta^*(\alpha)$.  Then there exists $C_1>0$ such that, for any positive $f$ harmonic on a ball $B = B_r$, $0<r<r_0$, in $\Pi =  M \times X$,  
$$
\sup_{\frac{1}{2} B} f \leq C_1 \inf_{\frac{1}{2}  B} f.
$$
\et

For our application below, we need to single out a particular result that is traditionally required in the course of proving this theorem. Fortunately, the standard proof of this result follows the Moser Iteration scheme, which poses no new difficulties in our mixed setting (that is, the computations are abstract and carry over without any new arguments). Thus we omit the details but send the interested reader to \cite{SC1}, Theorem 2.2.3, for the continuous case and to \cite{Delmotte1}, Proposition 5.3, for the discrete case. See also the paper \cite{CoulGrigor1} of Coulhon and Grigor'yan.
\bp \label{sobolev}
For every $\delta \in (0,1)$, there is a constant $C_2 = C_2(\delta)>0$ such that for every positive $f$ harmonic on $B = B_r$, $0<r<r_0$, in $M \times X$,
$$
\sup_{\delta B} | f | \leq C_2 \left ( \pi(B)^{-1} \int_{B} f^2 d\pi \right )^{1/2}.
$$
\ep
The proof in \cite{SC1} derives this result with $2$ replaced by $p \in (0,2]$, where $C_2$ must then depend on $p$. Thus, to complete the proof of the elliptic Harnack inequality, we seek a bound for the $L^p$ norm of the harmonic function $f$ over $B$ by $\inf_B f$. Although deriving the bound of Proposition \ref{sobolev} for $\sup f$ poses no new difficulty in the mixed setting, this $L^p$ estimate requires careful reasoning, especially with respect to the graph component. Much of the issues one faces are sourced in executing the chain rule in the discrete calculus. Indeed, Thierry Delmotte discusses in \cite{Delmotte2, Delmotte1} that the property $\Delta^*(\alpha)$ is the key assumption in deriving inequalities where one typically applies the chain rule in the original continuous setting (for example, Cacciopoli-type inequalities).

Given the linear nature of the product operator $\Delta = \Delta_M+\Delta_X$, the reader may believe that the computational subtleties used to derive the Elliptic version on the graph may carry over to the product space $\Pi = M \times X$. Unfortunately, these arguments, which are given in \cite{Delmotte1}, will not work here because they would critically rely upon the uniform bounds
$$
\frac{1}{C} \leq \frac{f(x,u)}{f(x,v)} \leq C
$$
for $u \sim v$ and for some $C>0$. The property $\Delta^*(\alpha)$ ensures these bounds hold in the purely discrete case, but the mixed space fails to inherit such strong bounds. Fortunately, arguments that Thierry Delmotte provides in \cite{Delmotte3} for the stronger parabolic Harnack inequality work well for the product structure. 

Because most steps in the proof of the elliptic case rely on abstract arguments whose details can be found in \cite{Delmotte1, SC1} and other sources, the authors believe it sufficient to exemplify how the arguments of Delmotte in \cite{Delmotte3} transfer to the product case. To this end, we have chosen to cover the details of a crucial step in proving that if $f$ is positive and harmonic on a ball $B$, then its $L^p$ average over $\delta B$ can be compared with $\inf_{\delta B} f$ (see Theorem 2.3.1 of \cite{SC1} for the statement). An abstract result of Bombieri-DeGuisti (Lemma 2.2.6 of \cite{SC1}) greatly simplified the original proof of this result. One of the (two) main assumptions of this lemma involves bounding the growth of $\log f$ relative to its mean. To achieve this, one usually applies the Poincar\'{e} inequality to $\log f$ and must bound the resulting $\nabla (\log f)$ term. Normally, a standard argument involving integration by parts and Cauchy-Schwartz completes the bound, but this argument relies upon applying the chain rule to $\log (f)$. In the discrete elliptic case, Delmotte uses the uniform bounds discussed above to circumvent this chain rule issue (see the proof of Lemma 3.2 in \cite{Delmotte1}); however, as already mentioned, these bounds do not necessarily hold in our mixed setting. Therefore we borrow an idea of Delmotte's argument in \cite{Delmotte3} to show how the desired bound may still be derived.

\bp \label{assump}
Let $B'' = \frac{1}{4} B$. Then there exists $C_3>0$ such that, for any positive $f$ harmonic on $B$ and any $\lambda > 0$,
$$
\pi( B'' \cap \{|\log f - c_3|  > \lambda \}) \leq C_3 \frac{\pi(B)}{\lambda}.
$$
where $c_3 = (\log f)_{B''}$, the average of $\log f$ over $B''$.
\ep

\bpf

Fix $(z,w) \in M \times X$. Suppose $f$ is a positive harmonic function on $B = B((z,w),r)$ and write $r'=\frac{1}{2} r$. Choose a test function $\psi(x,u)$ as follows. Let $\psi$ have support in $B' = \frac{1}{2} B$, $0 \leq \psi \leq 1$, and $|\nabla \psi| \leq 1/r'$. Further suppose $|\psi(x,u) - \psi(x,v)| \leq 1/r'$ if $u \sim v$ and $|\psi(x,u)| \leq  1/r'$ if $ \sum_{v \sim u : v \notin B_X} \mu_{uv} \neq 0$. These properties are satisfied by the choice
$$
\psi(x,u) = \left ( 1 - \frac{d((z,w),(x,u))}{r} \right)^+.
$$
Let $\phi(x,u) = \psi^2(x,u)/f(x,u)$.

We begin our work on a product $B_M \times B_X$, where $B_M = B(x,r/\sqrt{2})$ and $B_X = B(u,r/\sqrt{2})$. Further let $B'_M \times B'_X := \frac{1}{\sqrt{2}} B_M \times \frac{1}{\sqrt{2}}  B_X$ and notice that 
$$
B' \subset B'_M \times B'_X \subset B_M \times B_X \subset B
$$ 
By Definition \ref{harmonic} of harmonic function, $f$ satisfies
$$
\ba
0 & = \int_{B_M} \int_{B_X} \nabla_M \phi(x,u) \cdot \nabla_M f(x,u) + \int_{B_M} \int_{B_X} \phi(x,u) (- \Delta_X f(x,u)) = \mathcal{C}+\mathcal{D},
\ea
$$
where $\mathcal{C}$ and $\mathcal{D}$ denote the terms involving continuous and discrete differential operators, respectively. For the first term $\mathcal{C}$, we compute
$$
\ba
 \mathcal{C} & = \int_{B_M} \int_{B_X} \nabla_M \phi \cdot \nabla_M f  \\
&= \int_{B_M} \int_{B_X}  (2 \psi \nabla_M \psi \cdot \nabla_M (\log f)) - \psi^2 | \nabla_M(\log f) |^2)
\ea
$$
With the first term, we use Cauchy-Schwarz and the inequality $ab \leq \frac{1}{2}(\delta^{-1} a^2+\delta b^2)$, for some small $\delta>0$, to get 
\begin{equation}
\begin{split} \nonumber
& \left | \int_{B_M} \int_{B_X} 2 \nabla_M \psi \cdot \nabla_M (\log f)  \psi \right | \\
& \leq  \left ( \int_{B_M} \int_{B_X}   4 | \nabla_M \psi |^2 \right)^{\frac{1}{2}} \left (  \int_{B_M} \int_{B_X}   \psi^2 | \nabla_M (\log f) |^2 \right )^{\frac{1}{2}}\\
& \leq 4 \delta^{-1} \left ( \int_{B_M} \int_{B_X}    | \nabla_M \psi |^2 \right) + \delta \left (  \int_{B_M} \int_{B_X}  \psi^2 | \nabla_M (\log f) |^2 \right )
\end{split}
\end{equation}
These last two observations together give 
\begin{equation} \label{eqn1}
\begin{split} 
(1-\delta) \int_{B_M} \int_{B_X} \psi^2 |\nabla_M(\log f)|^2 \leq 4 \delta^{-1} \left ( \int_{B_M} \int_{B_X}   | \nabla_M \psi |^2 \right)
\end{split}
\end{equation}
Now for the discrete part $\mathcal{D}$, we must take the integration by parts carefully:
\begin{equation} \label{eqn2}
\begin{split} 
 \sum_{u \in B_X}  \phi(x,u) & \left ( \sum_{v \sim u} \mu_{uv} ( f(x,u) - f(x,v))\right )\\ 
 & =  \sum_{u \in B_X} \sum_{v \sim u: v \in B_X} \phi(x,u) \mu_{uv} ( f(x,u) - f(x,v)) \\ 
 & \ \ \ \ +  \sum_{u \in B_X} \sum_{v \sim u : v \notin B_X} \mu_{uv} \phi(x,u) f(x,u)  \\ 
 & \ \ \ \ -  \sum_{u \in B_X} \sum_{v \sim u : v \notin B_X} \mu_{uv} \phi(x,u) f(x,v)  \\
 & \leq   \sum_{u \in B_X} \sum_{v \sim u: v \in B_X} \phi(x,u) \mu_{uv} ( f(x,u) - f(x,v)) \\ 
 & \ \ \ \ + \sum_{u \in B_X} \phi(x,u) f(x,u) \sum_{v \sim u : v \notin B_X} \mu_{uv} 
\end{split}
\end{equation}
Next we exploit the (crucial) inequality (2.14) of Delmotte in \cite{Delmotte3}: 
\begin{equation} \nonumber
\begin{split} 
& \left (  \frac{\psi^2(x,u)}{f(x,u)}  -  \frac{\psi^2(x,v)}{f(x,v)} \right )(f(x,u) - f(x,v))  \leq  \\ 
& \left (36(\psi(x,u)-\psi(x,v))^2 - \frac{1}{2} \min \{\psi^2(x,u), \psi^2(x,v) \} \frac{(f(x,u)-f(x,v))^2}{f(x,u)f(x,v)} \right )
\end{split}
\end{equation}
This in turn gives us the following estimate on the discrete part:
\begin{equation} \nonumber
\begin{split} 
& \int_{B_M} \sum_{u \in B_X} \phi(x,u) \left ( \sum_{v \in B_X} \mu_{uv} ( f(x,u) - f(x,v) \right ) \\
& \ \ \ + \int_{B_M} \sum_{u \in B_X} \psi^2(x,u)  \sum_{v \sim u : v \notin B_X} \mu_{uv}  \\
&= \int_{B_M} \frac{1}{2} \sum_{u \in B_X} \sum_{v \in B_X} \mu_{uv} \left ( \frac{\psi^2(x,u)}{f(x,u)} - \frac{\psi^2(x,v)}{f(x,v)} \right)(f(x,u) - f(x,v)) \\
& \ \ \ +\int_{B_M} \sum_{u \in B_X} \psi^2(x,u) \sum_{v \sim u : v \notin B_X} \mu_{uv} 
\end{split}
\end{equation}
\begin{equation} \label{eqn3}
\begin{split} 
& \leq \int_{B_M} \frac{1}{2} \sum_{u \in B_X} \sum_{v \in B_X} \mu_{uv} \bigg ( 36(\psi(x,u)-\psi(x,v))^2  \\ &   \ \ \ - \int_{B_M} \frac{1}{2} \min \{\psi^2(x,u) , \psi^2(x,v) \} \frac{(f(x,u)-f(x,v))^2}{f(x,u)f(x,v)} \bigg ) \\
& \ \ \ +\int_{B_M}  \sum_{u \in B_X} \psi^2(x,u) \sum_{v \sim u : v \notin B_X} \mu_{uv}.
\end{split}
\end{equation}
Recall that $B'_M \times B'_X \subset B_M \times B_X $ and note $\psi^2(x,u) \geq ((2-\sqrt{2})/2)^2 \geq 1/16$ for $(x,u) \in B'_M \times B'_X$. Putting \eqref{eqn1}, \eqref{eqn2}, and \eqref{eqn3} together, we get
\begin{equation}
\begin{split} \nonumber
& (1-\delta) \int_{B'_M} \int_{B'_X} \psi^2 |\nabla_M(\log f)|^2 +  \int_{B'_M}   \frac{1}{2^6} \sum_{u \in B'_X} \sum_{v \in B'_X} \mu_{uv} \frac{(f(x,u)-f(x,v))^2}{f(x,u)f(x,v)} \\
& \leq 4 \delta^{-1} \left ( \int_{B_M} \int_{B_X}   | \nabla_M \psi |^2 \right)  +18 \int_{B_M}  \sum_{u \in B_X} \sum_{v \in B_X} \mu_{uv} (\psi(x,u)-\psi(x,v))^2 \\
&\ \ \ \ \ \ + \int_{B_M}  \sum_{u \in B_X} \psi^2(x,u) \sum_{v \sim u : v \notin B_X} \mu_{uv}.
\end{split}
\end{equation}
Using Calculus, one can check that $(\log x)^2 \leq \frac{(x-1)^2}{x}$ and deduce that
$$
(\log f(x,u) - \log f(x,v) )^2 \leq \frac{(f(x,u)-f(x,v))^2}{f(x,u)f(x,v)}.
$$
Recall the conditions on our test function $\psi$: $|\nabla(\psi)| \leq 1/r'$, $|\psi(x,u) - \psi(x,v)| \leq 1/r'$ if $u \sim v$ and $|\psi(x,u)| \leq  1/r'$ if $ \sum_{v \sim u : v \notin B_X} \mu_{uv} \neq 0$. Further let $m = \min \{  1-\delta, 1/2^6 \} $. All these facts together then give

\begin{equation}
\begin{split} \nonumber
m \int_{B'} & |\nabla(\log f)|^2 \leq m \int_{B'_M} \int_{B'_X} |\nabla(\log f)|^2  \\
& \leq  4 \delta^{-1} \int_{B_M} \int_{B_X} | \nabla_M \psi |^2 +18 \int_{B_M}  \sum_{u \in B_X} \sum_{v \in B_X} \mu_{uv} (\psi(x,u)-\psi(x,v))^2\\
&+ \int_{B_M}  \sum_{u \in B_X} \psi^2(x,u) \sum_{v \sim u : v \notin B_X} \mu_{uv} \\
 & \leq  \frac{C_3}{(r')^2}\pi (B) \\
\end{split}
\end{equation}
for some $C_3>0$. Now we know by our work in the previous section that our product space $\Pi = M \times X$ satisfies the (weak) Poincar\'{e} Inequality. Recall $B'' = \frac{1}{2} B' = \frac{1}{4} B$ and write $c_3 = (\log f)_{B''}$ for the average of $\log f$ over $B''$. Then there is some $P>0$ such that
$$	
\int_{B''} |\log f -c_3 |^2 d \pi \leq P (r')^2 \int_{B'} |\nabla (\log f ) |^2 d \pi
$$
At last, we may conclude there is a constant $C_3>0$ such that
$$
\lambda \pi(B'' \cap \{|\log f - c_3 | \geq \lambda \}) \leq C_3 \pi(B),
$$
which completes the proof.

\epf
As explained above, the proof technique for Proposition \ref{assump} together with classical arguments prove the main result Theorem \ref{elliptic} of this section.

\section{Kernel estimates on graphs}
\label{Kernel Estimates on Graphs}

Let $X$ be a countable graph with measure $m(u) = \sum_v \mu_{uv}$, where $\mu_{uv} = \mu_{vu} \geq 0 $ is a symmetric weight on $X \times X$. We assume that $(X,\mu_{uv})$ satisfies the geometric properties of volume doubling $V(D_X, r_0)$ and the weak Poincar\'{e} Inequality $P(P_X, r_0)$, for some $r_0 \in (0, \infty]$; see Section \ref{Geometric Properties} for definitions. Let $P(u,v) = \frac{\mu_{uv}}{m(u)}$ and define the discrete Markov Kernel recursively by $P^0(u,v) = \delta_{u}(v)$ and for $n \geq 1$,
$$
P^n(u,v) := \sum_w P(u,w) P^{n-1}(w,v).
$$
Though not symmetric, this kernel is easily shown to be reversible, i.e.,
$$
\frac{P^n(u,v)}{m(v)} = \frac{P^n(v,u)}{m(u)}.
$$
Normalizing, one may think of $p(u,v):=P(u,v)/m(v)$ as the discrete ``heat kernel".
We let this kernel act on functions in $l^2(X)$ by 
$$
Pf(u) := \sum_v p(u,v) f(v) m(v).
$$
Reversibility then translates to self-adjointness of this operator. Note that with this definition, $I-P = -\Delta_X$.

Next we introduce notation to discuss functions of the kernel operator $P$. Suppose that $F$ is real analytic at $1$ and write $F(1-t) = \sum_n a_n t^n$ with  $\sum |a_m|< \infty$. Let $K = K_F := F(I-P)$. Then we define the action of $K$ on functions in $l^2(X)$ as
$$
Kf(u) = \sum_v k(u,v) f(v) m(v)
$$
where
$$
k(u,v) = \sum_n a_n p^n(u,v).
$$
Similarly to the above, we write $K(u,v) = m(v) k(u,v)$. In the application below, we employ Carne's transmutation formula in \cite{Carne1}, which requires notation for a simple random walk on $\Z$. Following  \cite{PDSC1} and \cite{Carne1}, let $a \in [-1,1)$ and write $X^a_n$ for a simple random walk on $\Z$ with parameters
$$
\mathbf{P}(X^a_n =  \pm 1 |X_{n-1}) = \frac{1-a}{4}, \ \ \ \mathbf{P}(X_n^a = 0 | X_{n-1}^a)=\frac{1+a}{2}.
$$
Let $F : \R \to \R$ be a real analytic function as above, so that $F(1-t) = \sum_k a_k t^k$ with $\sum |a_m|< \infty$. Let $\beta$ be the Bernoulli measure (i.e. $\beta(\pm 1) = 1/2$) and write
$$
\beta_s := (1-s) \delta_0 + s \beta,
$$
where $a \in (-1,1)$ and $2s = 1-a$. Moreover, write $\beta^{(n)}$ for the $n$th convolution power of $\beta$, with the convention that $\beta^{(0)} = \delta_0$. Now let $f_s$ be the convolution kernel of the operator $F(I-\beta_s)$:
$$
f_s = \sum_n a_n \beta_s^{(n)}.
$$
Define $\Delta_s f = f * (\delta_0 - \beta_s)$. Then $\forall j$, 
$$
f_{s,j} : = \Delta^j_s f_s = \sum_n a_n (\delta_0-\beta)^j * \beta_s^{(n)}.
$$
The next definition defines the class of functions to which we will apply our result.
\bd
Fix $\psi: \N \to \R_+$, $s \in (0,1]$. The class $\mathcal{F}(s, r,\psi)$ is the set of functions $F$ analytic at $1$ such that, for all $j, q \in \N$,
$$
 \sum_{|m| \geq q} |f_{s,j}(m)| \leq \left ( \frac{2j}{r}\right )^{2j} \psi(q).
$$
\ed
Compare the next proposition to (2.4) in \cite{PDSC1}.

\bp \label{estimate}
Let $(X, \mu_{uv})$ be a symmetrically weighted graph satisfying $V(D_X, r_0)$. Fix $s \in (0,1]$ and $r$ with $0<r<r_0$. Assume that $\sigma(P) \subset (a,1], a \in (-1,1)$ with $1-a \leq 2s$. Suppose $F \in \mathcal{F}(s, r, \psi)$. Let $w_1,w_2 \in X$ satisfy $\rho:= d(w_1,w_2)=2r+q$, for $ q \in \N$. Let $B_1:=B_r(w_1), B_2 := B_r(w_2)$. Then we have 
$$
| \langle (I-P)^k K_F (I-P)^l \phi_1, \phi_2 \rangle | \leq \left ( \frac{2k+2l}{r}\right )^{2k+2l} \psi(q) \Vert \phi_1 \Vert_2 \Vert \phi_2 \Vert_2.
$$
for all $\phi_1, \phi_2$ with support in $B_1, B_2$, respectively.  
\ep
\bpf
As in Carne's paper \cite{Carne1}, introduce the Chebyshev polynomials
$$
Q_m(z): = \frac{1}{2}( (z+\sqrt{z^2-1})^m+ (z-\sqrt{z^2-1})^m)
$$
$$
Q_{a,m}(z) := Q_m \left ( \frac{2z-1-a}{1-a} \right ).
$$
Note that a change of variables $ z \to (w+w^{-1})/2$ gives
$$
Q_m(z) = (w^m+w^{-m})/2,
$$
which implies that $Q_m$ is bounded by $1$ on $[-1,1]$. A standard argument then implies that $Q_m(P)$ is a contraction on $l^2(X, m)$ (see \cite{Carne1} for details). Moreover, this computation along with symmetry properties of the simple random walk above lead to the formula
$$
P^n = \sum_m \mathbf{P}_0(X^a_n = m) Q_{a,m}(P).
$$
Now write $F(1-s) = \sum_{i=0}^{\infty} a_i s^i$ and observe
$$
\ba
K_F & = \sum_n a_n P^n = \sum_m \left ( \sum_n a_n \mathbf{P}_0(X^a_n = m) \right )Q_{a,m}(P) \\
&=\sum_m f_s(m) Q_{a,m}(P)
\ea
$$
Moreover,
$$
(I-P)^k K_F (I-P)^l  = \sum_m f_{s, l+k } (m) Q_{a,m}(P).
$$
This immediately gives us
$$
\langle (I-P)^k K_F (I-P)^l  \phi_1, \phi_2 \rangle = \sum_m f_{s,l+k}(m) \langle Q_{a,m}(P) \phi_1, \phi_2 \rangle.
$$
But since $Q_{a,m}(P)$ is a contraction, $| \langle Q_{a,m}(P) \phi_1, \phi_2 \rangle| \leq \Vert \phi_1 \Vert_2 \Vert \phi_2 \Vert_2 $. In addition, our set up implies that $p^i(u,v) = 0$ if $d(u,v) > i$. Thus $ \langle Q_{a,m}(P) \phi_1, \phi_2 \rangle = 0$ for $|m| < q$. Hence, we may conclude
$$
|\langle (I-P)^k K_F (I-P)^l  \phi_1, \phi_2 \rangle | \leq \Vert \phi_1 \Vert_2 \Vert \phi_2 \Vert_2  \sum_{|m| \geq q} | f_{s,l+k}(m) |,
$$
which establishes the desired inequality.
\epf
Finally, the next theorem provides the main application of our mixed elliptic Harnack inequality. Our argument follows the paper \cite{CGT1} of Cheeger, Gromov, and Taylor. As in that paper, we follow the notation that $\Vert \cdot \Vert_E$ and $| \cdot |_E$ are the $L^2$ and $L^{\infty}$ norms over the set $E$, respectively.

\bt \label{kernel}
Assume $(X,\mu)$ satisfies $V(D_X,r_0)$ and $P(P_X, r_0)$, $r_0 \in (0,\infty]$. Retain the assumptions and notation as in Proposition \ref{estimate}. Then there is a $C_4$ such that
$$
|k|_{B_{r/2}(w_1) \times B_{r/2}(w_2)} \leq  \frac{C_4 }{ \sqrt{m(B_{r}(w_1)) m(B_{r}(w_2))}} \psi(q),
$$
for all $F \in \mathcal{F}(s,r,\psi)$ and $r<r_0$.
\et

\bpf
We begin by collecting a few estimates. Let $F \in \mathcal{F}(s,r,\psi)$ and write $K = K_F$. Recall $d(w_1,w_2) = 2r+q$. Let $\phi \in l^2(X)$ and $\text{supp} \  \phi \subset B_r(w_2)$. For $|x| \leq r/2$ and $u \in B_r(w_2)$, define
$$
\xi(x,u) := \sum_{k=0}^{\infty} \frac{x^{2k}}{(2k)!} (I-P)^{k} (K \phi) (u),
$$
which is in $L^2([-r/2,r/2] \times B_r(w_1))$. The reader may recognize the (formal) identity $\xi(x,u) = \cosh[ x (\sqrt{ I-P})] (K \phi) (u)$, and note that 
$$
\left ( \frac{\partial^2}{\partial x^2} +  \Delta_{X} \right )\xi(x,u) = 0.
$$
In particular, $\xi(0,u) =  K \phi(u) = F(I-P) \phi(u)$. Proposition \ref{estimate} tells us that
$$
\Vert (I-P)^k K \phi \Vert_{B_r(w_1)} \leq  \left ( \frac{2k}{r} \right)^{2k} \psi(q) \Vert \phi \Vert_{B_r(w_2)}.
$$
This in turn implies there is some $C>0$ such that
$$
\Vert \xi  \Vert_{[-r,r] \times B_{r}(w_1)} \leq C \psi(q) \Vert \phi \Vert_{B_r(w_2)} r^{1/2}
$$
(here and below, the norm applies to suppressed arguments). This estimate along with Proposition \ref{sobolev} applied to the harmonic function $\xi$ together say there is a $C'>0$ such that
$$
\ba
|K\phi|_{B_{r/2}(w_1)} & = |\xi(0,\cdot )|_{B_{r/2}(w_1)}  \\
& \leq C' m(B_r(w_1))^{-1/2} r^{-1/2} \Vert \xi \Vert_{[-r/2,r/2] \times B_{r}(w_1)} \\ 
& \leq C C' m(B_r(w_1))^{-1/2} \psi(q) \Vert \phi \Vert_{B_{r}(w_1)}
\ea
$$
Maximizing over all such $\phi$ with $\Vert \phi \Vert = 1$, we have for each $u \in X$
\begin{equation} \label{est1}
\Vert k(u, \cdot ) \Vert_{B_{r}(w_2)} \leq C C' m(B_r(w_1))^{-1/2} \psi(q).
\end{equation}
Indeed, this same argument gives us, for each $u \in X$,
\begin{equation} \label{est2}
\Vert (I-P)^k k(u, \cdot ) \Vert_{B_{r}(w_2)} \leq C C' m(B_r(w_1))^{-1/2} \left ( \frac{2k}{r}  \right )^{2k} \psi(q).
\end{equation}
Similarly to the above, for each $u \in X$, define a function by
$$
\zeta_u(x,v):=\sum_{k=0}^{\infty}  \frac{x^{2k}}{(2k)!} (I-P)^k (K \delta_v)(u)
$$
and notice again that
$$
\left ( \frac{\partial^2}{\partial x^2} + \Delta_{X}  \right )\zeta_u(x,v) = 0.
$$
Noting $\zeta_u(0,v) = k(u,v)$, we conclude again by Proposition \ref{sobolev} that there is a $C''>0$ such that, for any $u \in X$,
\begin{equation} \label{est3}
\begin{split}
|k(u,\cdot )|_{B_{r/2}(w_2)} & = |\zeta_u(0,\cdot ))|_{B_{r/2}(w_2)} \\ 
& \leq C'' m(B_r(w_2))^{-1/2} r^{-1/2} \Vert \zeta_u\Vert_{[-r,r] \times B_{r}(w_1)}.
\end{split}
\end{equation}
Now we can bound the term $\Vert \zeta_u \Vert_{[-r/2,r/2] \times B_{r}(w_1)}$ with estimates \eqref{est1}, \eqref{est2} to get
\begin{equation} \label{est4}
\Vert \zeta_u \Vert_{[-r,r] \times B_{r}(w_1)}  \leq C C' m(B_{r}(w_1))^{-1/2} \psi(q) r^{1/2}.
\end{equation}
With the above results collected, we complete the proof directly. Combining the two estimates  \eqref{est3} and \eqref{est4}, we have, for each $u \in B_{r/2}(w_1)$, 
$$
\ba
|k(u,\cdot )|_{B_{r/2}(w_2)} & = |\zeta_u(0,\cdot )|_{B_{r/2}(w_2)}  \\
& \leq C'' m(B_r(w_2))^{-1/2} r^{-1/2} \Vert \zeta_u\Vert_{[-r,r] \times B_{r}(w_1)} \\
& \leq C'' m(B_r(w_2))^{-1/2} r^{-1/2} \cdot C C' m(B_{r}(w_1))^{-1/2} \psi(q) r^{1/2} \\
&  = C C' C'' [m(B_{r}(w_1)) m(B_{r}(w_2))]^{-1/2} \psi(q).
\ea
$$
Since this holds for all $u \in B_{r/2}(w_1)$, we conclude there is a constant $C_4>0$ such that
$$
|k|_{B_{r/2}(w_1) \times B_{r/2}(w_2)} \leq \frac{C_4}{ \sqrt{m(B_{r}(w_1)) m(B_{r}(w_2))}}\psi(q)
$$
\epf
Let us look at a concrete example to get a better feel for the statement of Theorem  \ref{kernel}. This example uses some results on convolution powers of functions on $\Z$ from the paper \cite{PDSC1} by Diaconis and Saloff-Coste.

\begin{example} \label{example}

Fix $n,k, l \in \N$ and retain the notation of Proposition \ref{estimate}. Assume that the weight $\mu$ induces a kernel $P$ satisfying $\sigma(P) \subset (a,1]$ and $s = (1-a)/2<2^{-1+1/k}$. Consider the function
$$
F(x) = x^l (1-x^k)^n.
$$
We show that $F$ is in $\mathcal{F}(s,r,\psi)$, where 
$$
\psi(q) = \left (\frac{1}{n^{l/k}} \right ) \exp \left (-c_5 \left ( \frac{q}{n^{1/2k}} \right )^{2k/(2k-1)}  \right )
$$
for some $c_5>0$ and where $r = O(n^{1/2k})$ (to be determined somewhat more precisely below).

Note $F(I-P) = (I-P)^l(I-(I-P)^k)^n$ and $$f_{s} = (\delta_0-\beta_s)^{*l}*(\delta_0-(\delta_0 - \beta_s)^{*k})^{*n}.$$ Let $g_{s} = \delta_0-(\delta_0 - \beta_s)^{*k}$ and notice $\widehat{g^{(n)}_s}(\theta) = (1-s^k(1-\cos \theta)^k)^n$. Moreover, given our choice of $s$, we have $|\widehat{g^{(n)}_s}(\theta)| < 1$. Lastly, it can be shown (see the discussion in Section 3 of \cite{PDSC1}) that
$$
\widehat{g^{(n)}_s}(\theta) = e^{ - n(s/2)^k \theta^{2k}(1+o(1))}.
$$
Observe that $g^{(n)}_s$ meets the technical assumptions of Theorem $3.3$ of \cite{PDSC1}, which tells us there exist $C_5, c_5 >0$ such that
$$
|\Delta^{j+l} g^{(n)}_s(m)| \leq \frac{C_5^{2(j+l)}}{n^{(j+l)/k}n^{1/2k}} \exp \left (-c_5 \left ( \frac{|m|}{n^{1/2k}} \right )^{2k/(2k-1)}  \right )
$$
Set $E_0(n) = [0, 2n^{1/2k})$ and $E_i(n) = [2^{i-1} n^{1/2k}, 2^{i} n^{1/2k})$ for $i \geq 1$. Then we can lastly estimate
\begin{equation}
\begin{split} \nonumber
\sum_{|m| \geq q} & \exp \left (-c_5 \left ( \frac{|m|}{n^{1/2k}} \right )^{2k/(2k-1)}  \right ) \\
& \leq  \exp \left (-\frac{c_5}{2} \left ( \frac{q}{n^{1/2k}} \right )^{2k/(2k-1)}  \right )  \sum_{|m| \geq q} \exp \left (-\frac{c_5}{2} \left ( \frac{|m|}{n^{1/2k}} \right )^{2k/(2k-1)}  \right ) \\
& \leq \exp \left (-\frac{c_5}{2} \left ( \frac{q}{n^{1/2k}} \right )^{2k/(2k-1)}  \right ) \sum_{i=0}^{\infty}  \sum_{m: |m| \in E_i(n) } \exp \left (-\frac{c_5}{2} \left ( 2^p \right )^{2k/(2k-1)}  \right ) \\
& \leq \exp \left (-\frac{c_5}{2} \left ( \frac{q}{n^{1/2k}} \right )^{2k/(2k-1)}  \right )  n^{1/2k} \sum_{p=0}^{\infty} 2^{p+1} \exp \left (-\frac{c_5}{2} \left ( 2^p \right )^{2k/(2k-1)}  \right ) \\
& \leq C_5  n^{1/2k}  \exp \left (-\frac{c_5}{2} \left ( \frac{q}{n^{1/2k}} \right )^{2k/(2k-1)}  \right ).
\end{split}
\end{equation}
Putting everything together, we can estimate the sum
$$
\ba
\sum_{|m| \geq q} |f_{s, l+k}(m)|  & = \sum_{|m| \geq q} |\Delta_s^{l+s} g^{(n)}_{s}(m)| \\ 
& \leq \sum_{|m| \geq q} \frac{C_5^{2(j+l)}}{n^{(j+l)/k}n^{1/2k}} \exp \left (-c_5 \left ( \frac{|m|}{n^{1/2k}} \right )^{2k/(2k-1)}  \right ) \\
& \leq  \frac{C_5^{2(j+l)}}{n^{(j+l)/k}n^{1/2k}} \cdot n^{1/2k} \cdot  \exp \left (-c_5 \left ( \frac{q}{n^{1/2k}} \right )^{2k/(2k-1)}  \right ) \\
& \leq \left( \frac{C_5}{n^{1/2k}} \right)^{2j+2l}    \exp \left (-c_5 \left ( \frac{q}{n^{1/2k}} \right )^{2k/(2k-1)}  \right ) \\
& \leq \left( \frac{2j}{r} \right)^{2j} \psi(q),
\ea
$$
where we have taken $r = n^{1/2k}/C_5^{2l}$. Hence, $F \in \mathcal{F}(s,r,\psi)$ and Theorem \ref{kernel} gives us the following kernel estimate:
$$
|k |_{B_{r/2}(w_1) \times B_{r/2}(w_2)} \leq \frac{C_5}{ n^{l/k} \sqrt{m(B_{r}(w_1)) m(B_{r}(w_2))}} \exp \left (-c_5 \left ( \frac{q}{n^{1/2k}} \right )^{2k/(2k-1)}  \right ).
$$
We summarize the result of this example in a slightly more explicit form.
\end{example}
\bt \label{explicit}
Fix $n,k,l \in \N$ and retain the same notation established in Proposition \ref{estimate}. Assume $r = n^{1/2k} < r_0$. Suppose $\sigma(P) \subset (a,1]$ and $(1-a)/2 < 2^{-1+1/k}$. Then there exist $C_5, c_5 >0$ such that
\begin{equation} \nonumber
\begin{split}
   \big |(I-P)^l & (I-(I-P)^k)^n(u,v) \big |   \\
& \leq \frac{C_5m(v) }{ n^{l/k} \sqrt{m(B_{r}(u)) m(B_{r}(v))}} \exp \left (-c_5 \left ( \frac{d_X(u,v)}{n^{1/2k}} \right )^{2k/(2k-1)}  \right ).
\end{split}
\end{equation}
\et
\bn
The case of most interest is when $r_0 = \infty$, i.e., the Volume Doubling property and weak Poincar\'{e} inequality hold globally. In this case, Theorem \ref{explicit} holds for all $r \in (0,\infty)$.
\en
\bn

Note that the addition of volume terms is the primary difference between this inequality and those found in the paper \cite{PDSC1} of Diaconis and Saloff-Coste. In contrast, the bound above \emph{does} capture the decay of $P^n_k(u,u)$, as discussed after Theorem 4.2 of \cite{PDSC1}.  

In addition, equation \eqref{vdeqn} of Lemma \ref{vdlemma} tells us that
$$
\frac{m(B_{r}(u))}{m(B_{r}(v))} \leq \frac{m(B_{r+d_X(u,v)}(v))}{m(B_{r}(v))} \leq D_X \left ( \frac{r+d_X(u,v)}{r} \right )^{\log D_X/\log 2}.
$$
These terms can be absorbed into the constant $c_5$ of the exponential, giving us the simplified bound
$$
\frac{C_5m(v) }{ n^{l/k} m(B_{r}(u))} \exp \left (-c_5 \left ( \frac{d_X(u,v)}{n^{1/2k}} \right )^{2k/(2k-1)}  \right )
$$
that involves only a single volume term.

\en

\section{Application to convolutions on groups} \label{applications}

The kernel estimates of the previous section have a number of interesting consequences for Markov Chains on groups. We need to fix some notation for discrete groups that will vary slightly from what was given above for graphs. Let $G$ be a group and $S$ a (finite) set of generators. Suppose further that $\mu$ is a probability measure supported on $S$ and that $\mu$ is symmetric in the sense that $\mu(u) = \mu(u^{-1}) =: \check{\mu}(u)$. The (left-invariant) random walk on $G$ driven by $\mu$ has transition kernel $P_{\mu}(u,v) = \mu(u^{-1}v)$. $P_{\mu}$ induces a kernel operator on $l^2(G)$ (with respect to counting measure) associated to convolution by $\mu$:
$$
P_{\mu}f(u) = \sum_{v \in G} f(v) \mu(v^{-1}u) = f*\check{\mu}(u).
$$
Observe
$$
\sum_v P_{\mu}(u,v) = \sum_v \mu(u^{-1}v) = \sum_v \mu(v) = 1.
$$
The iterated kernel is $P^n_{\mu}(u,v) = \mu^{(n)}(u^{-1}v)$ and it's invariant distribution can be given by $\pi \equiv 1/| \text{supp}(\mu)|$. As before, let $d_G$ denote the shortest path distance and write $|g| := d_G(e,g)$, where $e$ is the identity. Lastly, we let $V(n)$ denote the volume of an $n$-dimensional ball.

\begin{thm} \label{estimate2}

Let $G$ be a infinite discrete group and $S$ its set of generators. Endow $G$ with a probability measure $\mu$ and define a (symmetric) weight by $\mu_{uv} = \mu(vu^{-1})$. Let $\mu$ be finitely supported and generating. Suppose also that the induced reversible Markov kernel $P_{\mu}$ has spectrum of the form $(a,1] \subset [-1,1]$ with $(1-a)/2 < 2^{-1+1/k}$. Further assume that $G$ has polynomial volume growth, i.e., $V(n) \asymp n^d$. Define 
$$
\mu_k := \delta_e - (\delta_e - \mu)^{*k}.
$$
Then there exist constants $C_6, c_6 >0$ such that
$$
|\mu_k^{(n)}(g) | \leq \frac{C_6}{V(n^{1/2k})} \exp \left ( - c_6 \left ( \frac{|g|}{n^{1/2k}} \right )^{2k/(2k-1)} \right ).
$$
\end{thm} 
\bpf
This follows from our work in Example \ref{example} by taking $l=0$, i.e., by using the bound for the kernel $K_F$ where $F(x) = (1-x^k)^n$. Then there exist constants $C_6, c_6 >0$ 
$$
|\mu_k^{(n)}(g) | = |k(e,g)| \leq \frac{C_6}{V(n^{1/2k})} \exp \left ( - c_6 \left ( \frac{|g|}{n^{1/2k}} \right )^{2k/(2k-1)} \right ),
$$
which is our desired result.
\epf
\bc \label{cor1}
For any $k$ such that $(1-a)/2<2^{-1+1/k}$, there is a constant $N_k \in (0, \infty)$ such that, for any $n \in \N$,
$$
\sum_{g \in G} | \mu_k^{(n)} (g) | < N_k.
$$
\ec
This result is particularly interesting since, a priori, we know that $s = \sum_{g \in G} | \mu_k (g) | >1$, so a rough estimate would give
$$
\sum_{g \in G} | \mu_k^{(n)} (g) | \leq s^n,
$$
which goes to infinity as $n \to \infty$.
\bpf[Proof of Corollary \ref{cor1}]
As in Example \ref{example}, set $E_0(n) = [0, 2 n^{1/2k})$ and $E_i(n) = [2^{i-1} n^{1/2k}, 2^{i} n^{1/2k})$ for $i \geq 1$. Using the estimate of Theorem \ref{estimate2} and the volume growth assumption, we sum over all $g \in G$ to get
$$
\ba
\sum_{g \in G} | \mu_k^{(n)} (g) | &  \leq \sum_{g \in G}  \frac{C_6}{V(n^{1/2k})} \exp \left ( - c_6 \left ( \frac{|g|}{n^{1/2k}} \right )^{2k/(2k-1)} \right ) \\
& \leq \sum_{i=0}^{\infty} \sum_{ g :  |g| \in E_i(n)}   \frac{C_6}{V(n^{1/2k})} \exp \left ( - c_6 \left ( 2^m \right )^{2k/(2k-1)} \right ) \\
& \leq \sum_{m=0}^{\infty} C_6 \  \frac{V(2^{m+1} n^{1/2k})}{V(n^{1/2k})} \exp \left ( - c_6 \left ( 2^m \right )^{2k/(2k-1)} \right ) \\
& \leq \sum_{m=0}^{\infty} C_6 \  2^{(m+1)d} \exp \left ( - c_6 \left ( 2^{m} \right )^{2k/(2k-1)} \right ) < N_k,
\ea
$$
where the volume comparison term is implicitly absorbed into the constant $C_6$.
\epf

\bc
For any $k$ such that $(1-a)/2<2^{-1+1/k}$, there exists $C_6, c_6>0$ such that for all $n \in \N$ and $g \in G$,
$$
|\mu_k^{(n)}(g)-\mu_k^{(n+1)}(g)| \leq \frac{1}{n}  \frac{C_6}{V(n^{1/2k})} \exp \left ( - c_6 \left ( \frac{|g|}{n^{1/2k}} \right )^{2k/(2k-1)} \right ).
$$
Moreover, there is a constant $N_k \in (0, \infty)$ such that, for any $n \in \N$,
$$
\sum_{g \in G} | \mu_k^{(n)}(g) - \mu_k^{(n)}(g)| < N_k.
$$
\ec
\bpf
This again follows from our work in Example \ref{example}, but now with the choice $l = k$. To see this, observe
$$
\ba
\mu_k^{(n)}(g)-\mu_k^{(n+1)}(g) & = (\delta_e-\mu_k)*\mu_k^{(n)}(g) \\
& = (\delta_e - \mu)^{*k}*(\delta_e - (\delta_e - \mu)^{*k})^{*n},
\ea
$$
which corresponds to the case $F(x) = x^k(1-x^k)^n$. The other bound follows in exactly the same manner as the analogous bound of Corollary \ref{cor1}.

\epf

{ }

\end{document}